\documentclass[12pt,leqno]{amsart}
\newtheorem{dummy}{}[section]
\pdfoutput=1
\newtheorem{definition}[dummy]{Definition}
\newtheorem{theorem}[dummy]{Theorem}

\newtheorem{proposition}[dummy]{Proposition}
\newtheorem{lemma}[dummy]{Lemma}

\newtheorem{example}[dummy]{Example}
\begin{document}
\bibliographystyle{plain}
\title{Galois descent of determinants in the ramified case }
\author{Victor P. Snaith}
\date{October 2010; Mathematics subject classification: 19B28, 11S23 }

\begin{abstract}
In the local, unramified case the determinantal functions associated to the group-ring of a finite group satisfy Galois descent. This note examines the obstructions to Galois determinantal descent in the ramified case.
\end{abstract}
\maketitle
\section{Introduction}

This note arose out of an attempt to answer a question posed to me by Otmar Vejakob in the autumn of 2009.
I refer the reader to \cite{IV10} for background literature and a description of the related context of non-commutative Iwasawa theory, which is currently a very active research topic in algebraic number theory.

Let $M$ be a $p$-adic local field and $O_{M}$ its valuation ring. Let $G$ be a finite group and $O_{M}[G]^{*}$ the units in its $O_{M}$ group-ring. The determinant gives a homomorphism into Galois equivariant, unit-valued functions on $R(G)$, the complex representation ring of $G$ (see \S2). The determinantal functions are those in the image of this homomorphism, ${\rm Det}(O_{M}[G]^{*})$.

If $M/{\mathbb Q}_{p}$ is an unramified extension then there is an isomorphism of the form
\[   {\rm Det}(O_{M}[G]^{*})^{{\rm Gal}(M/ {\mathbb Q}_{p})}  \cong  {\rm Det}({\mathbb Z}_{p}[G]^{*})\].
This Galois descent isomorpism was first proved by M.J. Taylor \cite{Fro} using the Oliver-Taylor logarithm.
Using Explicit Brauer Induction \cite{SnEBI} I gave a simpler construction of the group-ring logarithm which, in turn, simplified the proof of unramified Galois descent for determinantal functions.

Otmar Venjakob's question was whether determinantal Galois descent held when $M/K$ was ramified; a question motivated by the case of $p$-adic non-commutative Iwasawa algebras.

In the case when $M/K$ is unramified the determinantal Galois descent hinges on the determination of $\hat{\alpha}_{G}({\rm Det}( 1 + A_{M}(G))$ where $\hat{\alpha}_{G}$ is defined in Proposition \ref{Eg3}. This in turn hinges on an integrality result concerning the image of the group-ring logarithm. In Proposition \ref{det6} I generalise this integrality result to the case where $M/K$ is arbitrary. The result gives the first signs of the difficulties which obstruct determinantal Galois descent - difficulties which are manifested in \cite{IV10}. In the general case the logarithmic image considered in Proposition \ref{det6} involves in an essential wayan element $h_{M} \in O_{M}$, which may not even be ${\rm Gal}(M/K)$-fixed. When $M/{\mathbb Q}_{p}$ is unramified Galois descent depends crucially on the fact that $h_{M}=p$, which is fixed by the Galois action.

A second problem in the general case is that my version of the group-ring logarithm depends on a choice of lifted Frobenius. This dependence makes the equivariance of the group-ring logarithm more complicated (see Proposition \ref{det8}).

After all these extra complications in the general case, in \S3 I can only offer a very modest example of ramified determinantal Galois descent, which I have included to illustrate the intended structure of the determinantal Galois descent proof.

  Written in January 2010, this note was not needed in \cite{IV10} and accordingly I undertook to post it independently.

\section{ Determinantal Congruences }

Throughout this section let  $p$  be a prime and let  $G$ be a finite group of order  $n$. Let  $N/{\mathbb Q}_{p}$  be a finite Galois extension containing all the  $n$-th roots of unity and let $M/K$ be a Galois subextension. Let $O_{M}$ and $ \pi_{M} O_{M}$ denote the valuation ring of $M$ and its maximal ideal, respectively.  We shall consider the group of  Galois-equivariant, unit-valued functions               
${\rm Hom}_{{\rm Gal}(N/M)}( R(G) , O_{N}^{*} )$                           
where the complex representation ring $R(G)$  is identified with  $R_{N}(G) = K_{0}(N[G])$ and is therefore generated by representations of the form $   T : G \longrightarrow  GL_{u}(N)  $.               

We have a determinantal homomorphism               
 \[   Det: O_{M}[G]^{*} \longrightarrow                 {\rm Hom}_{{\rm Gal}(N/M)}( R(G) , O_{N}^{*} ) \]                            given by the formula  
\[  Det( \sum_{\gamma} \lambda_{\gamma} \gamma ) (T) =   det(                \sum_{\gamma} \lambda_{\gamma} T(\gamma ))   \in  O_{N}^{*} . \]                

Choose $F \in {\rm Gal}(M/K)$ which is a lift of the Frobenius automorphism of the residue fields. If  $ O_{K}/ \pi_{K} = {\mathbb F}_{p^{q}}$ then $F(z) \equiv z^{p^{q}}$ (modulo $\pi_{M}O_{M}$) for all $z \in  O_{M}$. If  $ \sum_{\gamma} \lambda_{\gamma} \gamma  \in  O_{M}[G]$  we may                therefore define              
\[   F( \sum_{\gamma} \lambda_{\gamma} \gamma ) = \sum_{\gamma}                F(\lambda_{\gamma}) \gamma    \in  O_{M}[G]  ,  \]                
so that $F$ is a ring automorphism of  $O_{M}[G] $. The following result generalises (\cite{SnEBI} Theorem 4.3.10).                                   
\begin{theorem}{$_{}$}                
\label{det1}   
\begin{em}   
 
         Let  $z \in  O_{M}[G]^{*}$. Then, for all $ T \in R(G)$,                 
\[    Det(F(z))( \psi^{p^{q}}(T))/ ( Det(z)(T))^{p^{q}}  \in  1  +  \pi_{M}O_{N} .  \]
 Here  $\psi^{m}$  denotes the $m$-th Adams operation on $R(G)$. 
\end{em}               
\end{theorem}
\vspace{2pt}                              

{\bf Proof}  
\vspace{2pt} 

 There exist integers, $\{ n_{i} \} $ and one-dimensional representations 
 \linebreak
$ \{ \phi_{i} : H_{i}  \longrightarrow N^{*} \}$ such that, in $R(G)$,                              
\[    T =  \sum_{i}  n_{i} Ind_{H_{i}}^{G}( \phi_{i} )    \    {\rm and}    \                 \psi^{p^{q}}(T) =  \sum_{i}  n_{i} Ind_{H_{i}}^{G}( \phi_{i}^{p^{q}} ) . \]
By multiplicativity  we have             
\[    \frac{Det(F(z))(\psi^{p^{q}}(T))}{(Det(z)(T))^{p^{q}}} =   \prod_{i}      \          \frac{Det(F(z))(Ind_{H_{i}}^{G}(\phi_{i}^{p^{q}}))^{n_{i}}}{(Det(z)(Ind_{H_{i}}^{G}(\phi_{i})))^{p^{q}n_{i}}}                  \in  O_{N}^{*}                \]                
so that we are reduced to the comparison (modulo $ \pi_{M}O_{N}$) of the expressions               
 \[              Det( \sum_{\gamma} \lambda_{\gamma} \gamma)( Ind_{H_{i}}^{G}( \phi_{i}))^{p^{q}}                  \hspace{5pt}  {\rm and}  \hspace{5pt}                Det( \sum_{\gamma} F(\lambda_{\gamma}) \gamma)( Ind_{H_{i}}^{G}( \phi_{i}^{p^{q}}))                  \]                
in  $ O_{N}^{*}$, where  $ z =   \sum_{\gamma} \lambda_{\gamma} \gamma  \in  O_{M}[G]^{*}$.  Let us abbreviate  $(H_{i},\phi_{i}) $ to  $(H,\phi) $. Choose coset representatives,  $x_{1}, \ldots , x_{d} \in G$, for  $G/H$. There is a homomorphism, $ \sigma : G \longrightarrow \Sigma_{d}$, such that for  $g \in G$                               $gx_{i}  =  x_{ \sigma(g)(i)} h(i,g) $ with $ h(i,g) \in H$. Realising these two representations on the vector space $N[G] \otimes_{N[H]} N$,  with this notation we find that   
\[ Ind_{H}^{G}(\phi)(z) =  \sum_{\gamma} \lambda_{\gamma} \sigma(\gamma) \cdot diag[                \phi(h(1, \gamma)), \ldots , \phi(h(d, \gamma)) ] ) = X \] 
where  $diag[ u_{1}, \ldots , u_{d}]$  is the diagonal matrix whose  $(i,i)$-th entry is equal to $u_{i}$. Also we have   
\[  Ind_{H}^{G}(\phi^{p^{q}})(F(z)) =   \sum_{\gamma} F(\lambda_{\gamma})                \sigma(\gamma) \cdot diag[ \phi^{p^{q}}(h(1, \gamma)), \ldots , \phi^{p^{q}}(h(d,                \gamma)) ] )  = Y.  \]
The $(u,v)$-th entries of $X$ and $Y$ are given by
\[  X_{u,v} = \sum_{\gamma,  \   \sigma(\gamma)_{u,v}=1} \  \lambda_{\gamma} \phi(h(v, \gamma))  \  {\rm and}  \  Y_{u,v} = \sum_{\gamma,  \   \sigma(\gamma)_{u,v}=1} \   F(\lambda_{\gamma}) \phi(h(v, \gamma))^{p^{q}}  \]
respectively. Therefore
\[ \begin{array}{llr}
X_{u,v}^{p^{q}} &   \equiv  \sum_{\gamma,  \   \sigma(\gamma)_{u,v}=1} \  \lambda_{\gamma}^{p^{q}} \phi(h(v, \gamma))^{p^{q}}   &  ({\rm modulo }  \   pO_{N})  \\
\\
&  \equiv  \sum_{\gamma,  \   \sigma(\gamma)_{u,v}=1} \   F(\lambda_{\gamma}) \phi(h(v, \gamma))^{p^{q}}  &  ({\rm modulo }  \   \pi_{M} O_{N})  \\
\\
& =  Y_{u,v}  &
\end{array} \]
and so
\[   \begin{array}{llr}                det(X)^{p^{q}}    &    =     ( \sum_{ \beta \in \Sigma_{d}}  \   {\rm sign}(\beta)                X_{1,\beta(1)}^{l} \ldots  X_{d,\beta(d)})^{p^{q}}  &   \\
\\
&  \equiv   \sum_{ \beta \in \Sigma_{d}}  \   {\rm sign}(\beta)                X_{1,\beta(1)}^{p^{q}} \ldots  X_{d,\beta(d)}^{p^{q}}  & ({\rm modulo}  \    pO_{N})         \\       \\                   &  \equiv  \sum_{ \beta \in \Sigma_{d}}  \   {\rm sign}(\beta) Y_{1,\beta(1)} \ldots                 Y_{d,\beta(d)}   &  ({\rm modulo} \   \pi_{M} O_{N})  \\   
\\                &  =   det(Y)  &  ({\rm modulo}  \  \pi_{M}O_{N})                  \end{array}  \]      
which implies the result,  since  $det(X)$ and $ det(Y) $ both lie in the units of $O_{N}$. $\Box$
\begin{definition}            
\label{det2}                
\begin{em}   

In the situation of Theorem \ref{det1} we may define a homomorphism, which depends on the choice of Frobenius lift $F$, 
\[   {\rm Log}_{F}  :   O_{M}[G]^{*}  \longrightarrow  {\rm Hom}_{{\rm Gal}(N/M)}( R(G) ,  \pi_{M} O_{N})  \]
by the formula
\[    {\rm Log}_{F} (z)(T) =  log(Det(F(z))( \psi^{p^{q}}(T))/ ( Det(z)(T))^{p^{q}}) \]
where, for $x \in \pi_{M} O_{N}$, $log( 1 + x ) = \sum_{n=1}^{\infty} (-1)^{n-1} x^{n}/n \in  \pi_{M} O_{N}$ as usual. We shall also denote by $ {\rm Log}_{F} $ the composition
\[   {\rm Log}_{F}  :   O_{M}[G]^{*}  \longrightarrow  {\rm Hom}_{{\rm Gal}(N/M)}( R(G) ,  \pi_{M} O_{N})  \subset  {\rm Hom}_{{\rm Gal}(N/M)}( R(G) ,  N) . \]
             Denote by  $M \{ G \}$ the  $M$-vector space whose basis consists of the                conjugacy classes of elements of  $G$. Recall that there is an isomorphism of $M$-vector spaces (\cite{SnEBI} Proposition 4.5.14)
 \[  \psi : M\{G\} \stackrel{\cong}{\longrightarrow}  {\rm Hom}_{{\rm Gal}(N/M)}( R(G) ,  N)  \]   
 given by $\psi(\sum_{\gamma} \  m_{\gamma} \gamma)(T) =  \sum_{\gamma} \  m_{\gamma} {\rm Trace}(T(\gamma)) $.

The composition $  \psi^{-1} \cdot  {\rm Log}_{F} $ defines a logarithmic homomorphism, depending on the choice of Frobenius lift $F$, of the form
\[ \alpha_{G} :   O_{M}[G]^{*}  \longrightarrow M \{ G \} . \]

               The Jacobson radical, $ J \subseteq O_{M}[G]$   (\cite{Lang84} p.636) is the left ideal which is equal to the intersection of all the maximal left ideals of  $O_{M}[G]$. In fact, $J$, is a two-sided ideal and  $O_{M}[G]/J$  is semi-simple. Hence, by Wedderburn's theorem  (\cite{Lang84} p.629)  $O_{M}[G]/J$  is a product of matrix rings over division rings. However, in this finitely generated, local situation some power of  $J$  lies in  $ pO_{M}[G]$  (\cite{CR} vol. I, \S5.22, p.112). Hence the division algebras have characteristic  $p$.  Since the Brauer group of a finite field vanishes, there is an isomorphism of the form $O_{M}[G]/J $ is isomorphic to a finite product of rings of matrices with entries in finite fields of characteristic $p$.    

If $r \in J$ then $ r^{t} \in pO_{M}[G]$  for some positive integer $t$  (\cite{CR} vol.I,  \S5.22, p. 112)                and the series for  $ (1 - r)^{-1}$  converges  $p$-adically so that  $ 1 - r \in O_{M}[G]^{*}$.  Therefore $log(Det(1 - r)(T))$ and $log(Det(F(1-r))( \psi^{p^{q}}(T))$ both converge $p$-adically in $O_{N}$ and (c.f. \cite{SnEBI} Lemma 4.3.21)
\[  {\rm Log}_{F} (1-r)(T) =  log(Det(F(1-r))( \psi^{p^{q}}(T))) - log( ( Det(1-r)(T))^{p^{q}})). \]
Define $ L_{F,0}(1 - r) \in M[G]$ by the $p$-adically convergent series             \[ L_{F, 0}(1 - r) =  p^{q} \sum_{n=1}^{\infty} r^{n}/n  -                 \sum_{n=1}^{\infty} \hat{F}(r^{n})/n     \]                where  $ \hat{F}( \sum_{\gamma \in G} \lambda_{\gamma} \gamma ) =  \sum_{\gamma                \in G} F(\lambda_{\gamma}) \gamma^{p^{q}} $.      
\end{em}
\end{definition} 
\begin{proposition}{$_{}$}                
\label{det3}   
 \begin{em} 
             If $c : M[G] \longrightarrow  M\{G\}$ sends $\gamma \in G$ to its conjugacy class then, in the notation of Definition \ref{det2}, 
 \[           \alpha_{G}(1 - r) =  c(L_{F, 0}(1 - r))  \in M\{G\}. \]
 \end{em}               
\end{proposition}   
 \vspace{2pt}
 
 {\bf Proof}
 \vspace{2pt}             
             
Let  $ r \in J$  and  $ T \in R(G)$  be as in Definition \ref{det2}. Suppose that  $T$  is a representation and let  $ \lambda_{1}, \ldots , \lambda_{u}$  denote the eigenvalues of  $T(r)$. Each  $\lambda_{i}$  lies in the maximal ideal 
of $ O_{N}$ (\cite{SnEBI} Lemma 4.3.21) and therefore the following series converges:                             
\[  \begin{array}{ll}                  log((Det( 1 - T(r))^{p^{q}}))  &  =  p^{q}  log(Det( 1 - T(r)))  \\
\\
&=  p^{q}  log( \prod_{i=1}^{u} ( 1 - \lambda_{i} ))  \\                
 \\  &  =  -  p^{q}  \sum_{i=1}^{u} \sum_{m=1}^{\infty}  \lambda_{i}^{m}/m         \\        \\                 &  =  -  p^{q}   \sum_{m=1}^{\infty}  \sum_{i=1}^{u} \lambda_{i}^{m}/m        \\        \\                  &  =  -  p^{q}  \sum_{m=1}^{\infty} \  {\rm Trace}(T(r^{m}))/m .                 \end{array}  \]                
Similarly, since the eigenvalues of  $ \psi^{p^{q}}(T)(r)$  are  $ \{                 \lambda_{i}^{p^{q}} \}$, we find that                 \begin{displaymath}                log(Det(1 - F(r))( \psi^{p^{q}}(T)))  =  -  \sum_{m=1}^{\infty} {\rm Trace}(T( \hat{F}(r^{m})))/m  .                 
\end{displaymath}                Therefore $  \psi( \alpha_{G}(1 - r)) =  \psi(c(L_{F, 0}(1 - r)) )$ which completes the proof, since $\psi$ is an isomorphism.            $\Box$     
                         
\begin{definition}{$_{}$}              
\label{det4}                  
\begin{em}   
              
Define an $O_{M}$-submodule $ \Lambda_{G} $ of $M\{G\}$ by     
\[                 \Lambda_{G} = O_{M}[G]/( \sum_{x,y \in G} O_{M}(xy - yx)) ,              \]
which embeds into $M\{G\}$ via $c$, the homomorphism of Proposition \ref{det3}. For $ r \in J $ the element $c(L_{F, 0}(1 - r)) $ lies in $\Lambda_{G}$.

Define $h_{M/K} \in M$ (up to multiplication by $O_{M}^{*}$) by the formula
\[   h_{M}O_{M}  =  O_{M} \pi_{M}  \bigcup  O_{M}\pi_{M}^{p}/p  \bigcup   O_{M} \pi_{M}^{p^{2}}/p^{2}  \bigcup   \ldots  \bigcup O_{M}  \pi_{M}^{p^{k}}/p^{k}  \bigcup  \ldots . \]
If the ramification index of $M/{\mathbb Q}_{p}$ is equal to $e$, so that $p O_{M} = \pi_{M}^{e} O_{M}$, then
$h_{M} =  \pi_{M}^{{\rm min}_{k \geq 0 }(p^{k} - ke)}$. For example, if$M/{\mathbb Q}_{p}$ is unramified then $h_{M} = p$.
\end{em}                
\end{definition}                
\begin{lemma}{$_{}$}
\label{det5}
\begin{em}

For $x \in O_{M}$ and $k \geq 0$
\[   F( x^{p^{k}}) \equiv  x^{p^{q+k}}   \  ({\rm modulo}  \  {\rm H}(\pi_{M}^{p^{k}}, p \pi_{M}^{p^{k-1}},
p^{2} \pi_{M}^{p^{k-2}}, \ldots , p^{k} \pi_{M}) O_{M})  \]
where, in the notation of Definition \ref{det4}, 
\[   {\rm H}(\pi_{M}^{p^{k}}, p \pi_{M}^{p^{k-1}},
p^{2} \pi_{M}^{p^{k-2}}, \ldots , p^{k} \pi_{M}) =  p^{k}  \pi_{M}^{{\rm min}_{0 \leq i \leq k }(p^{i} - ie)}.  \]
\end{em}
\end{lemma}
\vspace{2pt}

{\bf Proof}
\vspace{2pt}    

By definition of the lifted Frobenius $F$ this is true for $k=0$. By induction suppose that
$ F( x^{p^{k}}) =  x^{p^{q+k}}   +  {\rm H}(\pi_{M}^{p^{k}}, p \pi_{M}^{p^{k-1}},
p^{2} \pi_{M}^{p^{k-2}}, \ldots , p^{k} \pi_{M}) \cdot z$ for some $z \in O_{M}$. Denoting ${\rm H}(\pi_{M}^{p^{k}}, p \pi_{M}^{p^{k-1}}, p^{2} \pi_{M}^{p^{k-2}}, \ldots , p^{k} \pi_{M})$ by $\lambda $ we have
\[ F(x^{p^{k+1}})   =  (  x^{p^{q+k}}  + \lambda z)^{p} =  x^{p^{q+k+1}}  + \lambda^{p}z^{p} + \sum_{j=1}^{p-1} \   
\binom{p}{j}  \ x^{jp^{q+k}}  \lambda^{p-j} z^{p-j}  \]      
so that $ F(x^{p^{k+1}})  -  x^{p^{q+k+1}} $ is congruent to zero modulo 
\[  {\rm H}(\pi_{M}^{p^{k+1}}, p^{p} \pi_{M}^{p^{k}},
p^{2p} \pi_{M}^{p^{k-1}}, \ldots , p^{kp} \pi_{M}^{p}, 
p \pi_{M}^{p^{k}}, p^{2} \pi_{M}^{p^{k-1}},
p^{3} \pi_{M}^{p^{k-2}}, \ldots , p^{k+1} \pi_{M} )  \]
which is ${\rm H}(\pi_{M}^{p^{k+1}}, p \pi_{M}^{p^{k}}, p^{2} \pi_{M}^{p^{k-1}},
p^{3} \pi_{M}^{p^{k-2}}, \ldots ,  p^{k+1} \pi_{M})$, as required. $\Box$
\begin{proposition}{$_{}$}               
\label{det6}     
\begin{em}

           Let  $G$  be any finite group. If  $ r \in J$  then $c(L_{F, 0}(1 - r)) $ lies in $h_{M} \cdot \Lambda_{G}$.        \end{em}         
\end{proposition}     
\vspace{2pt}                           

{\bf Proof} 
\vspace{2pt}  
     
   Consider the series                
\[              L_{F, 0}(1 - r) =  p^{q} \sum_{m=1}^{\infty} \  r^{m}/m  -  \sum_{m=1}^{\infty} \   \hat{F}(r^{m})/m  \in  M[G].                 \]   
                      
If  $p^{q}$  does not divide  $m$  then  $ p^{q} r^{m}/m  \in  pO_{M}[G]$  and the  $c(p^{q} r^{m}/m)  \in  p \Lambda_{G}$. Now consider the remaining terms in the series  
\[    \sum_{m=1}^{\infty} p^{q} (  r^{p^{q}m} / p^{q} m)  -   \hat{F}(r^{m})/m                    =  \sum_{m=1}^{\infty} ( r^{p^{q}m}  -  \hat{F}(r^{m}))/m .             \]                
Suppose that  $ m = p^{s}u$  with  $HCF(u,p) = 1$  then we may set  $ t =  r^{u}$  so that                
\[       ( r^{p^{q}m}  -  \hat{F}(r^{m}))/m  =  (t^{p^{q+s}}  -  \hat{F}(t^{p^{s}}))/ p^{s}u .  \]              
Therefore it will suffice to show that, if  $ t \in J$ and $s \geq 0$,                 
\[  c( t^{p^{q+s}}  -  \hat{F}( t^{p^{s}}))  \in  c( h_{M} p^{s} O_{M}[G]) = h_{M} p^{s} \Lambda_{G} .    \]
    
If $G = \{ g_{1}, \ldots , g_{n} \}$ suppose that  $ t = \sum_{i=1}^{n} a_{i}g_{i}  $ so that
\[ t^{p^{q+s}} = \sum_{\underline{j}}  \     a_{j_{1}} \ldots a_{j_{p^{q+s}}} g_{j_{1}} \ldots                g_{j_{p^{q+s}}}        \]
where $\underline{j}$ ranges over all possible $p^{q+s}$-tuples. The cyclic group $C$ of order $p^{q+s}$ acts on the set of $\underline{j}$'s by cyclic permutation. The products  $ g_{1} \ldots g_{v}$  and   $ g_{2} \ldots                g_{v}g_{1}$  are conjugate in  $G$  so that each term in the subsum of terms which are cyclically conjugate to $\underline{j}$ will have the same image under $c$. Similarly
\[  \hat{F}( t^{p^{s}}) =  \sum_{\underline{k}}  \     F(a_{k_{1}} \ldots a_{k_{p^{s}}})
 g_{k_{1}} \ldots  g_{k_{p^{s}}}  g_{k_{1}} \ldots  g_{k_{p^{s}}}   \ldots  g_{k_{1}} \ldots g_{k_{p^{s}}}     \]
 where the product $ g_{k_{1}} \ldots                g_{k_{p^{s}}}   $ is repeated $p^{q}$ times.  
 
 Note that if $\underline{j} = ( k_{1}, \ldots , k_{p^{s}}, k_{1}, \ldots , k_{p^{s}}  ,  \ldots  k_{1}, \ldots , k_{p^{s}})$, as for example in the above expression for $ \hat{F}( t^{p^{s}})$, then $p^{q}$ divides the stabiliser order of $\underline{j}$ in $C$.

Suppose that the stabiliser of $\underline{j}$ in $C$ has order $p^{w}$ with $0 \leq w \leq q-1$ then the $C$-orbit of $\underline{j}$ has order at least $p^{s+1}$ and the subsum consisting of these terms has image under $c$ which lies in $p^{s+1} \Lambda_{G}$. There are no terms in $ \hat{F}( t^{p^{s}})$ corresponding to such a $\underline{j}$.

Now suppose  that the stabiliser of $\underline{j}$ in $C$ has order $p^{w}$ with $q \leq w \leq q+s$. In this case
$\underline{j} =  (k_{1}, \ldots , k_{p^{q+s-w}}, k_{1}, \ldots , k_{p^{q+s-w}}, \ldots , k_{1}, \ldots , k_{p^{q+s-w}})$ where $k_{1}, \ldots , k_{p^{q+s-w}}$ is repeated $p^{w}$ times. Associated to the term
\[  a_{\underline{j}}g_{\underline{j}} =  a_{k_{1}} \ldots  a_{k_{p^{q+s-w}}} 
\ldots , a_{k_{1}}, \ldots , a_{k_{p^{q+s-w}}} g_{k_{1}} \ldots  g_{k_{p^{q+s-w}}} \ldots g_{k_{1}} \ldots g_{k_{p^{q+s-w}}} \]
is the term
\[ b_{\underline{j}} =  -   F( a_{k_{1}} \ldots  a_{k_{p^{q+s-w}}} \ldots a_{k_{1}} \ldots  a_{k_{p^{q+s-w}}} ) g_{k_{1}} \ldots  g_{k_{p^{q+s-w}}}
g_{k_{1}} \ldots  g_{k_{p^{q+s-w}}}  \]
in which $ a_{k_{1}} \ldots  a_{k_{p^{q+s-w}}} $ and $g_{k_{1}} \ldots  g_{k_{p^{q+s-w}}} $ are repeated $p^{w-q}$ times.
The cyclic group $C'$ of order $p^{s}$ acts on the $p^{s}$-tuple
\[ (k_{1}, \ldots , k_{p^{q+s-w}}, k_{1}, \ldots , k_{p^{q+s-w}}, \ldots, k_{1}, \ldots , k_{p^{q+s-w}} )  \] 
with stabiliser of order $p^{w-q}$. Therefore the image under $c$ of the $C$-orbit of $a_{\underline{j}}g_{\underline{j}} $ is equal to $p^{q+s}/p^{w}$ copies of $c(a_{\underline{j}}g_{\underline{j}} )$ and the image under $c$ of the $C'$-orbit
of $b_{\underline{j}}$ is $p^{s}/p^{w-q}$ copies of $c(b_{\underline{j}})$. Therefore the image under $c$ of the sum over these two orbits, associated to $\underline{j}$, is equal to
\[  p^{q+s-w}  c( ( a_{k_{1}} \ldots  a_{k_{p^{q+s-w}}} )^{p^{w}} - F(( a_{k_{1}} \ldots  a_{k_{p^{q+s-w}}} )^{p^{w-q}} ) g_{k_{1}} \ldots  g_{k_{p^{q+s-w}}} g_{k_{1}} \ldots g_{k_{p^{q+s-w}}} \]
which, by Lemma \ref{det5}, lies in
\[ p^{q+s-w}  {\rm H}(\pi_{M}^{p^{w-q}}, p \pi_{M}^{p^{w-q-1}},
p^{2} \pi_{M}^{p^{w-q-2}}, \ldots , p^{w-q} \pi_{M})  \Lambda_{G}  \subseteq p^{s}h_{M} \Lambda_{G}  \]
as required. $\Box$
\begin{dummy}{The ${\rm Gal}(M/K)$-equivariance}
\label{det7}
\begin{em}

Let $\sigma \in {\rm Gal}(M/K)$. Then $\sigma$ acts on $z = \sum_{\gamma} \ \lambda_{\gamma} \gamma \in  O_{M}[G]$ by the formula $\sigma(z) =  \sum_{\gamma} \   \sigma(\lambda_{\gamma} ) \gamma$. Also $\sigma$ acts on
$f \in {\rm Hom}_{{\rm Gal}(N/M)}( R(G) ,  N)$ by the formula $(\sigma \cdot f)(T) =  \tilde{\sigma}(f( \tilde{\sigma}^{-1}(T)))$, where $ \tilde{\sigma} \in {\rm Gal}(N/K)$ is any lifting of $\sigma$. If $F \in {\rm Gal}(M/K)$ is one choice of lifted Frobenius, as in Definition \ref{det2} then $\sigma F \sigma^{-1}$ is another.
The following result describes the relation between the homomorphisms $ {\rm Log}_{F}$ and $ {\rm Log}_{\sigma F \sigma^{-1}}$.
\end{em}                
\end{dummy}   
\begin{proposition}{$_{}$}                
\label{det8}     
\begin{em}
           In the notation of \S\ref{det7} $Log_{\sigma F \sigma^{-1}} ( \sigma(z)) = \sigma \cdot (Log_{F}(z))$.     \end{em}        
 \end{proposition}  
\vspace{2pt}                               

{\bf Proof} 
\vspace{2pt}   
           
If  $ \tilde{\sigma} \in {\rm Gal}(N/K)$ is a lift of $\sigma$, as in \S\ref{det7}, we have
\[   \begin{array}{l}
   \sigma \cdot ( Log_{F})(z))(T)  \\
   \\
   =    \tilde{\sigma} (   log(Det(F(z))( \psi^{p^{q}}(  \tilde{\sigma}^{-1}(T)))/ ( Det(z)( \tilde{\sigma} ^{-1}(T)))^{p^{q}})   )  \\
   \\
   =  log(  \frac{Det(  \sum_{\gamma} \   \tilde{\sigma}(F(\lambda_{\gamma}))   \tilde{\sigma}( \tilde{\sigma}^{-1}( T(\gamma^{p^{q}} ))   )       )}{Det(  \sum_{\gamma} \     \tilde{\sigma}(\lambda_{\gamma}) 
    \tilde{\sigma}( \tilde{\sigma}^{-1}(T(\gamma) ))  )^{p^{q}} }    )   \\
    \\
    =     log(  \frac{Det(  \sum_{\gamma} \   \tilde{\sigma}(F(\lambda_{\gamma}))   T(\gamma^{p^{q}}    )       )}{Det(  \sum_{\gamma} \     \tilde{\sigma}(\lambda_{\gamma}) 
  T(\gamma)   )^{p^{q}} }   )   \\
  \\
      =     log(  \frac{Det(  \sum_{\gamma} \   \tilde{\sigma}(F( \lambda_{\gamma}))   T(\gamma^{p^{q}}    )       )}{Det(  \sum_{\gamma} \     \tilde{\sigma}(\lambda_{\gamma}) 
  T(\gamma)   )^{p^{q}} }   )   \\
  \\
  = Log_{\sigma F \sigma^{-1}} ( \sigma(z))(T)  ,
   \end{array}  \]
as required. $\Box$


\section{A $2$-group example}
\begin{definition}                
\label{Eg1}                
\begin{em}   

             In this section we shall suppose that  we are in the situation of \S1 and that $G$  is a finite  $p$-group. Define  $A_{M}(G)$  to be equal to the kernel of the natural map from  $O_{M}[G]$  to  $O_{M}[G^{ab}]$
\[  A_{M}(G) = {\rm Ker}( O_{M}[G]  \longrightarrow  O_{M}[G^{ab}]) .  \]       
Therefore  $ A_{M}(G)$ is contained in $J$, the Jacobson radical of  $O_{M}[G]$, which was introduced in  \ref{det2}. We are going to study the subgroup  $ Det( 1 + A_{M}(G) )$  of  $ {\rm Hom}_{{\rm Gal}(N/M)}( R(G), O_{N}^{*} )$.               
\end{em}                
\end{definition}

Firstly we recall a result of  \cite{Wa79}.                
\begin{proposition}{$_{}$}                
\label{Eg2} 
\begin{em}  
               $Det( 1 + A_{M}(G) )$   is torsion free subgroup of 
\linebreak
${\rm Hom}_{{\rm Gal}(N/M)}( R(G) , O_{N}^{*} )$.    
\end{em}                
\end{proposition} 

\vspace{10pt}
   \begin{proposition}{$_{}$}                \label{Eg3}   
\begin{em}  
                        There is a well-defined, injective homomorphism    
\[                 \hat{\alpha}_{G} :   Det( 1 + A_{M}(G))  \longrightarrow M\{ G \}                \] 
given by the formula $ \hat{\alpha}_{G}(Det(u)) = \alpha_{G}(u)$, where $\alpha_{G}$ is the homomorphism introduced in  Definition \ref{det2}.
\end{em}                   
\end{proposition}  
\vspace{2pt}                              

{\bf Proof} 
\vspace{2pt} 

This result was proved in (\cite{SnEBI} \S4.5) when $M/ {\mathbb Q}_{p}$ was unramified, using the fact that $\alpha_{G}$ was ${\rm Gal}(M/ {\mathbb Q}_{p})$-equivariant. In the current, more general, situation this is no longer true because $M/K$ need not be linearly disjoint from the cyclotomic $p$-power extension $K(\mu_{p^{\infty}})/K$. Accordingly, we shall use instead the partial Galois equivariance of Proposition \ref{det8}.
                 When  $G$  is abelian there is nothing to prove; therefore we will assume that $G$  is non-abelian.           
                       Let  $ u \in 1 + A_{M}(G)  \subset O_{M}[G]^{*}$  and suppose that  $ \alpha_{G}(u) =                0$.  By Proposition \ref{det8}, for all $j \geq 0$,
\[     0 =  F^{j} \cdot  Log_{F}(u)  =  Log_{F}( F^{j}(u))  . \]
This occurs if and only if the homomorphism                
\[                 T \mapsto log( \frac{ Det(F^{j+1}(u)) \psi^{p^{q}}(T)}{ (Det(F^{j}u)(T))^{p^{q}}})  =  Log_{F}(F^{j}(u))(T)                \]               
 is zero for all $ T \in R(G)$. Therefore, since  $R(G)$  is finitely generated, there exists a positive integer $m$ such that for all $T$ and $j \geq 0$
\[            Det(F^{j+1}(u)) (\psi^{p^{q}}(T))^{p^{m}}   =  Det(F^{j}u)(T)^{p^{m + q }}   .             \]                  
Hence we have
\[ \begin{array}{l}
  Det(F^{j}(u)) (\psi^{jp^{q}}(T))^{p^{m}}    \\
  \\
  =    Det(F^{j-1}(u)) (\psi^{(j - 1)p^{q}}(T))^{p^{m +q}}    \\
  \\
  =    Det(F^{j-2}(u)) (\psi^{(j - 2)p^{q}}(T))^{p^{m + 2q}}    \\
  \\
  =    \hspace{40pt}   \vdots     \hspace{40pt}   \vdots     \hspace{40pt}   \vdots    \hspace{40pt}   \vdots  \\
  \\
  =    Det(u)(T)^{p^{m + jq}} 
\end{array}  \]                
                Now suppose that  $ \#(G) = j p^{q}$ with ${\rm HCF}(j,p)=1$ so that $ \psi^{j p^{q}}(T)  =  dim(T)  \in  R(G)$. This means that, if  $ \epsilon: O_{M}[G]  \longrightarrow  O_{M}[\{ 1 \}] =  O_{M} $  is the augmentation map, then                 
\[             Det( u^{p^{m + jq}})(T)    =  Det( \epsilon( F^{j}(u)))^{dim(T)}    =  1     ,           \]                 
since  $ \epsilon( 1 + A_{M}(G))  = \{ 1 \} $. Therefore  $ Det(u)^{^{p^{m + jq}}}  =  1 $  and so  $ Det(u)  =  1 $, by  Prposition \ref{Eg2},  which shows that $\hat{\alpha}_{G}$ is injective, provided that it is well-defined.      
                
 To show well-definedness suppose that $Det(u)(T) = Det(T(u))  = 1$ for all $T \in R(G)$. If $u = \sum \ m_{\gamma} \gamma$ then $F(u) =    \sum \ F(m_{\gamma}) \gamma$ and $T( F(u)) =      \sum \ F(m_{\gamma}) T(\gamma)$. Let $\tilde{F} \in {\rm Gal}(N/K)$ be a lift of $F$ and suppose that $T = \tilde{F}(T')$ then
 \[  Det(T(F(u)) = Det(     \sum \ F(m_{\gamma}) \tilde{F}(T')(\gamma)   )  = \tilde{F}(Det(T'(u))) = 1  .  \]
 Therefore $Log_{F}(u) = 0$ and $\alpha_{G}(u) = 0$. If $Det(u') = Det(u'')$ then 
 \linebreak
 $Det(u'(u'')^{-1}) = 1$ and so 
 $0 = \alpha_{G}(u'(u'')^{-1}) = \alpha_{G}(u') - \alpha_{G}(u'')$ and $\hat{\alpha}_{G}$ is well-defined.
    $\Box$    
 \vspace{15pt}

With the results obtained so far we do not have a complete description of  the image of $ Det( 1 + A_{M}(G))$ under $\alpha_{G}$, except in the unramified case, because Proposition \ref{det6} is more complicated in the presence of ramification. However, our results are sufficient for a very unambitious example of Galois descent for determinantal functions. 
\begin{example}{An easy case when $p=2$}
\label{Eg4}
\begin{em}

Let $G$ be a finite $2$-group containing a central element $z$ of order $2$ which is a commutator and such that 
$H = G/\langle z \rangle$ is abelian. In the notation of Definition \ref{det2}, in particular $p=2$ and the residue degree equals $2^{q}$, we shall prove next that     
\[    c( L_{F, 0}( 1 - (1 - z)O_{M}[G]))  =  2^{q} \cdot c( (1 - z)O_{M}[G]) .   \]    
If  $ x \in O_{M}[G]$  then, since  $z$  is central, $ \hat{F}((1 - z)x) = 0 $  because  $ z^{2} = 1 $. 
 Therefore                 
\[  \begin{array}{ll}                    L_{F, 0}( 1 - (1 - z)x )  &  =  2^{q} ( \sum_{n=1}^{\infty} (1 - z)^{n} x^{n}/n )      \\
 \\                                &    \equiv  2^{q} (1 - z)(x + x^{2})   \  ({\rm  modulo}  \   2^{q}(1 - z)^{2}x^{2})                 \end{array}  \]   
 so that  $ L_{F, 0}(1 - (1 - z)x )  \in  2^{q}(1 - z)O_{M}[G]  $ and therefore  
 \[  c(L_{F, 0}(1 - (1 - z))O_{M}[G])  \subseteq       2^{q}c( (1 - z))O_{M}[G]) .     \]
 If  $x \in J$  then one sees that  $   2^{q} c((1 - z)x)  \in  c(L_{F, 0}(1 - (1 - z)O_{M}[G])) $, by means of a standard approximation argument. This implies that       
\[    2^{q} c((1 - z)J)  \subseteq  c(L_{F, 0}(1 - (1 - z)O_{M}[G]))  \subseteq 
2^{q} \cdot c( (1 - z)O_{M}[G]) . \]
    
On the other hand we claim that               
\[   c( (1 - z)O_{M}[G])  \subseteq c((1 - z)J)    \]                
which will complete the proof. Write  $ z = a^{-1}b^{-1}ab$  for  $ a , b \in G$  and let  $ v = \sum m_{\gamma} \gamma  \in  O_{M}[G]$. We have to show that $c( (1 - z)v )   \in   c((1 - z)J) .$                Rewrite  $v$  as  $ v = \sum m_{\gamma} ( \gamma - a ) + \sum m_{\gamma} a$  so that                
\[  \begin{array}{ll}                  c( ( 1 - z )v )  &  =  c( (1-z)( \sum m_{\gamma} ( \gamma - a ) + \sum                m_{\gamma} a  ) )  \\  
\\                      & = c( (1-z) \sum m_{\gamma} ( \gamma - a ) )  ,                \end{array}  \]                
since  $  c( a - za ) = c( a - b^{-1}ab ) = 0 $.  However,  $ \gamma - a  \in                 Ker(O_{M}[G]  \longrightarrow  O_{M} ),$ the kernel of the augmentation,  and since  $G$  is an  $2$-group the kernel of the augmentation lies inside  $J$.               

Since, in this example, $ (1 - z)O_{M}[G] = A_{M}(G)$, we have established the following result.
\end{em}
\end{example}
\begin{proposition}{$_{}$}
\label{Eg5}
\begin{em}

In Example \ref{Eg4}                
\[    \hat{\alpha}_{G}( Det(1  +    A_{M}(G)))  =  2^{q} \cdot c( A_{M}(G)).  \]                   
\end{em}
\end{proposition}
\vspace{10pt}

We conclude this section by proving Galois descent for determinantal functions in the situation of Example \ref{Eg4}, following the argument given in (\cite{SnEBI} \S4.5) for the unramified case. 
 \begin{proposition}{$_{}$}                
\label{Eg6} 
 \begin{em}

                Let  $M/K$ be as in \S1 with $p=2$ and let $G$ be as in Example \ref{Eg4}. Then                
\[               Det(O_{M}[G]^{*})^{{\rm Gal}(M/K)}  \cong  Det( O_{K}[G]^{*}) .             \]
 \end{em}            
 \end{proposition} 
\vspace{2pt}

 {\bf Proof}  
\vspace{2pt}

For completeness, although it follows (\cite{SnEBI} \S4.5), we shall give the complete proof.             Consider the following diagram, whose rows are easily seen to be short exact.                \newline                                \begin{picture}(150,180)(-20,-50)                \put(-10,120){$ 1 + A_{M}(G)$}                \put(85,124){\vector(1,0){7}}                \put(110,120){$O_{M}[G]^{*}$}                \put(170,124){\vector(1,0){7}}                \put(190,120){$O_{M}[G^{ab}]^{*}$}                \put(10,90){\vector(0,-1){30}}                \put(120,90){\vector(0,-1){30}}                \put(200,90){\vector(0,-1){30}}                \put(20,75){Det}                \put(125,75){Det}                \put(205,75){Det}                \put(190,75){$\cong$}                \put(-10,30){$ Det(1 + A_{M}(G))$}                \put(85,34){\vector(1,0){7}}                \put(97,30){$Det(O_{M}[G]^{*})$}                \put(170,34){\vector(1,0){7}}                \put(185,30){$Det(O_{M}[G^{ab}]^{*})$}                \put(10,10){\vector(0,-1){15}}                \put(120,10){\vector(0,-1){15}}                \put(200,10){\vector(0,-1){15}}                \put(7,-25){$\{ 1 \}$}                \put(112,-25){$\{ 1 \}$}                \put(192,-25){$\{ 1 \}$}                \end{picture}                \newline                in which the vertical maps are induced by the determinant, which is an                isomorphism for abelian groups.                                Let  $U = G(M/K)$ then we may compare the bottom row for $K$ with                the  $U$-invariants of the bottom row for  $M$.                \newline                       \begin{picture}(150,150)(-5,0)                \put(-10,30){$ Det(1 + A_{M}(G))^{U}$}                \put(85,34){\vector(1,0){7}}                \put(100,30){$ Det(O_{M}[G]^{*})^{U}$}                \put(190,34){\vector(1,0){7}}                \put(205,30){$ Det(O_{M}[G^{ab}]^{*})^{U}$}                \put(10,90){\vector(0,-1){30}}                \put(120,90){\vector(0,-1){30}}                \put(220,90){\vector(0,-1){30}}                \put(20,75){$ \beta_{1}$}                \put(125,75){$ \beta_{2}$}                \put(225,75){$ \beta_{3}$}                \put(210,75){$\cong$}                \put(0,75){$\cong$}                \put(-10,120){$ Det(1 + A_{K}(G))$}                \put(85,124){\vector(1,0){7}}                \put(97,120){$Det( O_{K}[G]^{*})$}                \put(190,124){\vector(1,0){7}}                \put(205,120){$Det( O_{K}[G^{ab}]^{*})$}                \end{picture}                \newline                               The map, $ \beta_{3}$, is an isomorphism because                $ (O_{M}[G^{ab}]^{*})^{U}   \cong  O_{K}[G^{ab}]^{*}$                              and therefore the lower sequence is short exact. 

At this point the argument of (\cite{SnEBI} \S4.5) concludes by observing, since $\alpha_{G}$ and $\hat{\alpha}_{G}$ are 
Galois equivariant in that unramified case, that
$\beta_{1}$  may be identified with the isomorphism 
\begin{displaymath}                  2^{q} \cdot c( A_{K}(G))  \cong  2^{q} \cdot ( c( A_{M}(G)))^{U}                \end{displaymath}                so that  $ \beta_{2}$  is an isomorphism, by the five-lemma, which completes                the proof.   

However, in the general situation $\hat{\alpha}_{G}$ is {\em not} Galois equivariant in the naive sense; instead we must use Proposition \ref{det8}. Specifically, in Example \ref{Eg4}, we saw that $1  +   A_{M}(G) = 1 + (1-z) O_{M}[G]$ and that $L_{F,0}( 1 - (1-z)x)$ is independent of the choice of $F$ so that Proposition \ref{det8} (or the explicit formulae of Example \ref{Eg4}) imply that, in the present case, $\hat{\alpha}_{G}$ is Galois equivariant, as required. $\Box$



 \end{document}